\documentclass[a4paper,draft]{amsproc}
\usepackage{amssymb}

\theoremstyle{plain}

\theoremstyle{definition}

\numberwithin{equation}{section}
\renewcommand{\le}{\leqslant}
\renewcommand{\ge}{\geqslant}

\def\txt#1{{\textstyle{#1}}}
\def\hf{{\textstyle\frac{1}{2}}}
\def\a{\alpha}\def\b{\beta}
\def\d{{\,\rm d}}
\def\e{\varepsilon}
\def\f{\varphi}
\def\G{\Gamma}
\def\k{\kappa}
\def\s{\sigma}

\def\={\;=\;}

\def\zt{\zeta(\hf+it)}

\def\D{\Delta}

\def\R{\Re{\rm e}\,} 
\def\z{\zeta}

\def\txt#1{{\textstyle{#1}}}
\def\f{\varphi}
\font\tenmsb=msbm10
\font\sevenmsb=msbm7
\font\fivemsb=msbm5
\newfam\msbfam
\textfont\msbfam=\tenmsb
\scriptfont\msbfam=\sevenmsb
\scriptscriptfont\msbfam=\fivemsb
\def\Bbb#1{{\fam\msbfam #1}}

\font\teneufm=eufm10
\font\seveneufm=eufm7
\font\fiveeufm=eufm5
\newfam\eufmfam
\textfont\eufmfam=\teneufm
\scriptfont\eufmfam=\seveneufm
\scriptscriptfont\eufmfam=\fiveeufm
\def\mathfrak#1{{\fam\eufmfam\relax#1}}

\def \NN {\Bbb N}

\def \ZZ {\Bbb Z}

\title[Running title]{On the mean square of\\ 
the Riemann zeta-function in short intervals}

\subjclass[2000]{11 M 06; 11 N 37}


\author[Surname]{\bfseries Aleksandar Ivi\'c}

\address{
Katedra Matematike RGF-a\\
Universitet u Beogradu\\
\DJ u\v sina 7, 11000 Beograd\\
Serbia
}
\email{ivic@rgf.bg.ac.yu, aivic\_2000@yahoo.com}

\dedicatory{Communicated by }

\begin{document}

{\begin{flushleft}\baselineskip9pt\scriptsize
PUBLICATIONS DE L'INSTITUT MATH\'EMATIQUE\newline\strut\qquad\quad
Nouvelle s\'erie, tome ??(??) (2009), ??--??
\end{flushleft}}
\vspace{18mm}
\setcounter{page}{1}
\thispagestyle{empty}

\begin{abstract}
It is proved that, for $T^\e\le G = G(T) \le\hf\sqrt{T}$,
\begin{eqnarray*}
\int_T^{2T}\Bigl(I_1(t+G,G)-I_1(t,G)\Bigr)^2\d t &=& TG\sum_{j=0}^3a_j\log^j
\Bigl(\frac{\sqrt{T}}{G}\Bigr)\\&
+& O_\e(T^{1+\e}G^{1/2}+T^{1/2+\e}G^2)
\end{eqnarray*}
with some explicitly computable constants $a_j\;(a_3>0)$ where,
for fixed $k\in\Bbb N$,
$$
I_k(t,G) = {1\over\sqrt{\pi}}\int_{-\infty}^\infty |\z(\hf+it+iu)|^{2k}
{\rm e}^{-(u/G)^2}\d u.
$$
The generalizations to the mean square of $I_1(t+U,G) - I_1(t,G)$
over $[T,\,T+H]$  and the estimation of the mean square of $I_2(t+U,G)-I_2(t,G)$
are also discussed.
\end{abstract}

\maketitle

\section{Introduction}

The mean values of the Riemann zeta-function $\z(s)$, defined as
$$
\z(s) = \sum_{n=1}^\infty n^{-s}\qquad(\s = \R s > 1),
$$
(and otherwise by analytic continuation) occupy a central place in the
theory of $\z(s)$.
 Of particular significance is the mean square on the
``critical line" $\s = \hf$, and a vast literature exists on this
subject (see e.g., the monographs [3], [4], and [18]). One usually
defines the error term in the mean square formula for $|\zt|$ as
$$
E(T) := \int_0^T|\zt|^2\d t - T\Bigl(\log{T\over2\pi} + 2\gamma -1\Bigr),
\leqno(1.1)
$$
where $\gamma = -\G'(1) =0.57721 56649 \ldots$ is Euler's constant.
More generally, one hopes that for a fixed $k$ the function ($E(T) \equiv E_1(T)$
in this notation)
$$
E_k(T) := \int_0^T|\zt|^{2k}\d t - TP_{k^2}(\log T)\qquad(k\in \NN)\leqno(1.2)
$$
represents the error term in the asymptotic formula for the $2k$-th
moment of $|\zt|$, where $P_\ell(z)$
is a suitable polynomial in $z$ of degree $\ell$. This is  known, besides
the case $k=1$, only in the case $k=2$ (see e.g., [4], [17]), and any further
improvement would be of great significance, in view of numerous applications of power
moments of $|\zt|$. By means of random matrix theory plausible values of the
coefficients of the polynomial
$P_{k^2}(z)$ that ought to be standing in (1.2) are given by J.B. Conrey et al. [2].
However, these values are still conjectural.

\medskip
As for explicit formulas for $E_k(T)$ and related functions, we begin by mentioning
the famous formula of F.V. Atkinson [1] for $E(T)$.
Let $0 < A < A'$ be any two fixed constants
such that $AT < N < A'T$, let $d(n) =\sum_{\delta|n}1$ be the number of divisors
of $n$, and finally let $N' = N'(T) =
T/(2\pi) + N/2 - (N^2/4+ NT/(2\pi))^{1/2}$. Then
$$
E(T) = \sum\nolimits_1(T) + \sum\nolimits_2(T) + O(\log^2T),\leqno(1.3)
$$
where
$$
\sum\nolimits_1(T) = 2^{1/2}(T/(2\pi))^{1/4}\sum_{n\le N}(-1)^nd(n)n^{-3/4}
e(T,n)\cos(f(T,n)),\leqno(1.4)
$$
$$
\sum\nolimits_2(T) = -2\sum_{n\le N'}d(n)n^{-1/2}(\log (T/(2\pi n))^{-1}
\cos(T\log (T/(2\pi n)) - T + \pi /4),\leqno(1.5)
$$
with
$$\leqno(1.6)
$$
\begin{eqnarray*}
f(T,n)  &=& 2T{\rm {ar\,sinh}}\,\bigl(\sqrt{\pi n/(2T})\bigr) + \sqrt{2\pi nT
+ \pi^2n^2} - \txt{\frac{1}{4}}\pi\\
&=&  -\txt{1\over4}\pi + 2\sqrt{2\pi nT} +
\txt{1\over6}\sqrt{2\pi^3}\frac{n^{3/2}}{T^{1/2}} +
\b_5\frac{n^{5/2}}{T^{3/2}} +
\b_7\frac{n^{7/2}}{T^{5/2}} + \ldots\,,
\end{eqnarray*}
where the $
\b_j$'s are constants,
\[\leqno(1.7)\]
\begin{eqnarray*}
e(T,n) &=& (1+\pi n/(2T))^{-1/4}{\Bigl\{(2T/\pi n)^{1/2}
{\rm {ar\,sinh}}\,(\sqrt{\pi n/(2T})\,)\Bigr\}}^{-1}\\
&=& 1 + O(n/T)\quad\qquad(1 \le n < T),
\end{eqnarray*}
and $\,{\rm{ar\,sinh}}\,x = \log(x + \sqrt{1+x^2}\,).$
Atkinson's formula came into prominence several decades after its
appearance, and incited much research (see e.g., [3],[4] for some of them).
The presence of the function $d(n)$ in (1.4) and the structure
of the sum $\Sigma_1(T)$ point out the analogy between $E(T)$ and $\D(x)$,
the error term in the classical divisor problem, defined as
$$
\D(x) = \sum_{n\le x}d(n) - x(\log x + 2\gamma-1).
$$
This analogy was investigated by several authors, most notably by M. Jutila
[13], [14], and then later by the author [6]--[8] and [9]--[10].
Jutila [13] proved that
$$\leqno(1.8)$$
\begin{eqnarray*}
&{}& \int\limits_T^{T+H}\Bigl(\D(x+U)-\D(x)\Bigr)^2\d x \\
&=&{1\over4\pi^2}\sum_{n\le {T\over2U}}{d^2(n)\over
n^{3/2}}\int\limits_T^{T+H} x^{1/2}\left|\exp\left(2\pi iU\sqrt{{n\over
x}}\,\right)-1\right|^2\d x\\
& + &O_\e(T^{1+\e} +
HU^{1/2}T^\e),
\end{eqnarray*}
 for $1 \le U \ll T^{1/2} \ll H \le T$,
and an analogous result holds also for the integral of $E(x+U)-E(x)$
(the constants in front of the sum and in the exponential will be
$1/\sqrt{2\pi}$ and $\sqrt{2\pi}$, respectively). Here and later $\e\;(>0)$
denotes arbitrarily small constants, not necessarily the same ones at
each occurrence, while $a \ll_\e b$ means that the implied $\ll$--constant
depends on $\e$.
From (1.8) one deduces ($a\asymp b$ means that $a\ll b\ll a$)
$$
\int_T^{T+H}\Bigl(\D(x+U)-\D(x)\Bigr)^2\d x
\;\asymp\; HU\log^3\left({\sqrt{T}\over U}\right)
\leqno(1.9)
$$
for $HU \gg T^{1+\e}$ and $T^\e \ll U \le \hf\sqrt{T}$. In [14]
Jutila proved that the integral in (1.9) is
$$
\ll_\e T^\e(HU + T^{2/3}U^{4/3})\qquad(1 \ll H,U \ll T).
$$
This bound and (1.9) hold also for the integral of $E(x+U)-E(x)$.
Furthermore Jutila op. cit. conjectured that
$$
\int_T^{2T}\Bigl(E(t+U) - E(t-U)\Bigr)^4\d t \ll_\e T^{1+\e}U^2\leqno(1.10)
$$
holds for $1 \ll U \ll T^{1/2}$, and the analogous formula should hold for
$\D(t)$ as well. In fact, using the ideas of K.-M. Tsang [19] who
investigated the fourth moment of $\D(x)$, it can be shown that one
expects the integral in (1.10) to be of the order
$TU^2\log^6(\sqrt{T}/U)$. As shown in [11], the truth of Jutila's conjecture
(1.10) implies the hitherto unknown eighth moment bound
$$
\int_0^T|\zt|^8\d t \;\ll_\e\;T^{1+\e},\leqno(1.11)
$$
which would have important consequences in many problems from multiplicative
number theory, such as  bounds involving differences between consecutive primes.

\medskip
Despite several results on $E_2(T)$ (see e.g., [17]), no explicit
formula is known for this function,
which would be analogous to Atkinson's formula (1.3)-(1.7).
This is probably due to the complexity of the function in question,
and it is even not certain that such a formula exists. However,
when one works not directly with the moments of $|\zt|$, but with smoothed
versions thereof, the situation changes. Let, for $k\in\NN$ fixed,
$$
I_k(t,G) := \frac{1}{\sqrt{\pi}}\int_{-\infty}^\infty
|\z(\hf + it + iu)|^{2k}{\rm e}^{-(u/G)^2}\d u\qquad(1\ll G \ll t).
\leqno(1.12)
$$
Y. Motohashi's monograph [17] contains explicit formulas for $I_1(t,G)$
and $I_2(t,G)$ in suitable ranges for $G=G(t)$. The formula for $I_2(t,G)$
involves quantities from the spectral theory of the
non-Euclidean Laplacian (see op. cit.).
Let, as usual, $\,\{\lambda_j = \k_j^2 + {1\over4}\} \,\cup\,
\{0\}\,$ be the discrete spectrum of the non-Euclidean Laplacian
acting on $\,SL(2,\ZZ)\,$--automorphic forms, and $\a_j =
|\rho_j(1)|^2(\cosh\pi\k_j)^{-1}$, where $\rho_j(1)$ is the first
Fourier coefficient of the Maass wave form corresponding to the
eigenvalue $\lambda_j$ to which the Hecke series $H_j(s)$ is
attached. Then, for
$T^{1/2}\log^{-D}T \le G \le T/\log T$, and
for an arbitrary constant $D>0$, Motohashi's formula gives
$$
\leqno(1.13)
$$
\begin{eqnarray*}
&{}&G^{-1}I_2(T,G) = O(\log^{3D+9}T) \,+ \\&
+& \frac{\pi}{\sqrt{2T}}\sum_{j=1}^\infty \a_j H_j^3(\hf)\k_j^{-1/2}
\sin\left(\k_j\log{\k_j\over4{\rm e}T}\right)
\exp\Bigl(-{1\over4}\Bigl({G\k_j\over T}\Bigr)^2\Bigr).
\end{eqnarray*}
For our purposes the range for which (1.13) holds is not large enough.
We shall use a more precise the expression for $I_2(T,G)$ that
can be derived by following the proof of Theorem 3, given  in [12],
and then taking $\s\to\hf+0$. Namely if
$$
Y_0 = Y_0(T;\k_j):= {\k_j\over T}
\Biggl(\sqrt{1+{\Bigl({\k_j\over4T}\Bigr)}^2}
+ {\k_j\over2T}\Biggr),\leqno(1.14)
$$
then, for $\,T^\e \ll G= G(T) \ll T^{1-\e}\,$, it follows that
$$
 \leqno(1.15)
$$
\begin{eqnarray*}
&{}&I_2 (T,G)  \sim F_0(T,G) + O(1)
\\&
+ &{\pi G\over\sqrt{2T}}
\sum\limits_{\kappa_j
\le TG^{-1} \log T} \alpha_j  H^3_j(\hf)\kappa_j^{-1/2}
{\rm e}^{-{1\over4}G^2\log^2(1+Y_0)}\\&
\times& \sin \Bigl(\kappa_j \log
{\kappa_j\over 4{\rm e}T} +c_3\k_j^3T^{-2}+ \cdots + c_N\k_j^NT^{1-N}\Bigr).
\end{eqnarray*}
Here $N(\ge 3)$ is a sufficiently large integer, and all the constants
$c_j$ in (1.15) may be effectively evaluated.
The meaning of the
symbol $\sim$ is that besides the spectral sums in (1.15) a finite
number of other sums  appear, each of which is similar in nature
to the corresponding sum above, but of a lower order
of magnitude. The function $F_0(t,G)$ is given explicitly e.g., by eq. (5.112) of
[4]. We have
$$
F_0(t,G) = \R\left\{{G\over\sqrt{\pi}}\int_{-\infty}^\infty
\left[B_1{\G'\over\G}+\cdots +B_{11}{(\G')^2\G''\over\G^2}\right]
(\hf+it+iu)\,{\rm e}^{-(u/G)^2}\d u\right\},\leqno(1.16)
$$
where $B_1,\,\,\ldots\,,B_{11}$ are suitable constants.
The main contribution to $F_0(t,G)$ is of the form $GP_4(\log t)$, where
$P_4(z)$ is a polynomial of degree four in $z$, whose  coefficients
can be explicitly evaluated. This is easily
obtained by using the asymptotic formula
$$
{\G^{(k)}(s)\over\G(s)} =
\sum_{j=0}^kb_{j,k}(s)\log^js + c_{-1,k}s^{-1} +\ldots +c_{-r,k}s^{-r}
+ O_r(|s|^{-r-1})\leqno(1.17)
$$
for any fixed integers $k\ge1, r\ge0$, where each $b_{j,k}(s) \,(\sim
b_{j,k}$
for a suitable constant $b_{j,k}$) has an asymptotic expansion in
non-positive powers of $s$. One obtains (1.17)  from Stirling's
classical formula for $\G(s)$.

\section{Statement of results}

In [11] the author improved (1.9) and its analogue for $E(T)$ to a true
asymptotic formula. Namely it was shown that,
for $1 \!\ll\! U\! = \!U(T)\! \le \!\hf {\sqrt{T}}$, we have ($c_3 = 8\pi^{-2}$)
$$\leqno(2.1)$$
\begin{eqnarray*}
\int_T^{2T}\Bigl(\D(x+U)-\D(x)\Bigr)^2\d x & =& TU\sum_{j=0}^3c_j\log^j
\Bigl({\sqrt{T}\over U}\Bigr) \,+ \\
&+& O_\e(T^{1/2+\e}U^2) + O_\e(T^{1+\e}U^{1/2}),
\end{eqnarray*}
 a similar result being true
if $\D(x+U)-\D(x)$ is replaced by $E(x+U)-E(x)$, with different constants $c_j$.
It follows that, for $T^\e \le U = U(T) \le T^{1/2-\e}$,
(2.1) is a true asymptotic formula.
Moreover, for $T \le x \le 2T$ and $T^\e \le U= U(T) \le T^{1/2-\e}$,
from (2.1) it follows that
$$
\D(x+U) - \D(x) = \Omega \Bigl\{\sqrt{U}\log^{3/2}\Bigl({\sqrt{x}\over U}\Bigr)\Bigr\},\;
E(x+h) - E(x) = \Omega \Bigl\{\sqrt{U}\log^{3/2}\Bigl({\sqrt{x}\over U}\Bigr)\Bigr\}.
\leqno(2.2)
$$
These omega results ($f(x) = \Omega(g(x))$
means that $\lim_{x\to\infty}f(x)/g(x)\ne0$) show that Jutila's conjectures made in [13],
namely that
$$
\D(x+U) - \D(x) \ll_\e x^\e\sqrt{U},\;E(x+U) - E(x) \ll_\e x^\e\sqrt{U}\leqno(2.3)
$$
for $x^\e \le U \le x^{1/2-\e}$ are (if true), close to being best possible.
It should be also mentioned that the formula (2.1) can be further generalized, and the
representative cases are the classical circle problem and the summatory function of
coefficients of holomorphic cusp forms, which  were also treated in [11].

\medskip
In this work the problem of the mean square of the function
$I_k(t,G)$ (see (1.12)) in short intervals is considered
when $k = 1$ or $k=2$. The former case is much less
difficult, and in fact an asymptotic formula in the most important case
of the mean square of $I_1(t+G,G)-I_1(t,G)$ can be obtained.
The result, which is similar to (2.1), is

\bigskip
THEOREM 1. {\it For $\,T^\e\le G = G(T) \le\hf\sqrt{T}\,$ we have
$$\leqno(2.4)
$$
\begin{eqnarray*}&&
\int\limits_T^{2T}\bigl(I_1(t+G,G)-I_1(t,G)\bigr)^2\d t = \\
&=&TG\sum_{j=0}^3a_j\log^j
\Bigl({\sqrt{T}\over G}\Bigr)
+ O_\e(T^{1+\e}G^{1/2}+T^{1/2+\e}G^2)
\end{eqnarray*}
with some explicitly computable constants $a_j\;(a_3>0)$}.

\medskip
{\bf Corollary 1}. For $T^\e\le G = G(T) \le T^{1/2-\e}$, (2.4) is
a true asymptotic formula.

\medskip
{\bf Corollary 2}. For $T^\e\le G = G(T) \le T^{1/2-\e}$ we have
$$
I_1(T+G,G) - I_1(T,G) = \Omega(\sqrt{G}\log^{3/2}T). \leqno(2.5)
$$

\medskip
Namely if (2.5) did not hold, then replacing $T$ by $t$, squaring and integrating
we would obtain that the left-hand side of (2.4) is $o(TG\log^3T)$ as $T\to\infty$,
which contradicts the right-hand side of (2.4). The formula given by Theorem 1 makes
it then plausible to state the following conjecture, analogous to (2.3).

\medskip
{\bf Conjecture}. For $T^\e\le G = G(T) \le T^{1/2-\e}$ we have
$$
I_1(T+G,G) - I_1(T,G) = O_\e(T^\e\sqrt{G}).
$$

\medskip
The generalization of Theorem 1 to the mean square of $I_1(t+U,G) - I_1(t,G)$
over $[T,\,T+H]$ is will be discussed in Section 4, subject to the
condition (4.1). This is technically more involved than (2.5), so we
have chosen not to formulate our discussion as a theorem.

\medskip
The mean square of $I_2(t+U,G) - I_2(t,G)$ is naturally more involved than
the mean square of $I_1(t+U,G) - I_1(t,G)$. At the possible state of knowledge
involving estimates with $\k_j$ and related exponential integrals, it does not
seem possible to obtain an asymptotic formula, but only an upper bound. This is

\bigskip
THEOREM 2. {\it For} $T^\e \le U \le GT^{-\e}\ll T^{1/2-\e}, U = U(T), G=G(T)$
{\it we have}
$$
\int_T^{2T}\Bigl(I_2(t+U,G)-I_2(t,G)\Bigr)^2\d t \;\ll_\e \;
T^{2+\e}\left({U\over G}\right)^2.
\leqno(2.6)
$$

\section{The proof of Theorem 1 }

We shall prove first the bound
$$
\int_T^{2T}\Bigl(I_1(t+G)-I_1(t,G)\Bigr)^2\d t \ll TGL^3
\qquad(T^\e \le G \le \hf\sqrt{T}),
\leqno(3.1)
$$
where henceforth $L = \log T$ for brevity. This shows that the sum on the
left-hand side of (2.5) is indeed of the order given by the right-hand side,
a fact which will be needed a little later. We truncate (1.12) at $u =\pm GL$,
then differentiate (1.1) and integrate by parts to obtain
$$\leqno(3.2)$$
\begin{eqnarray*}
I_1(t,G) &=& {1\over\sqrt{\pi}}\int_{-GL}^{GL} \left(\log {t+u\over2\pi}
+ 2\gamma + E'(t+u)\right)
{\rm e}^{-(u/G)^2}\d u + O({\rm e}^{-L^2/2})\\
&=& {1\over\sqrt{\pi}}\int_{-GL}^{GL}
\left(\log {t+u\over2\pi} + 2\gamma + {2u\over G^2}E(t+u)\right)
{\rm e}^{-(u/G)^2}\d u + O({\rm e}^{-L^2/2}).
\end{eqnarray*}
This gives, for $T\le t \le 2T, T^\e\ll G\le \hf\sqrt{T}$,
$$
I_1(t+G,G)-I(t,G)= {2\over\sqrt{\pi}G^2}\int\limits_{-GL}^{GL}u{\rm e}^{-(u/G)^2}
\Bigl(E(t+u+G)-E(t+u)\Bigr)\d u + O\Bigl({G^2\over T}\Bigr).\leqno(3.3)
$$
We square (3.3), noting that $G^2/T \ll1$, and then integrate over $[T,\,2T]$.
We use the Cauchy-Schwarz inequality
for integrals together with Jutila's bound (1.9) (or (2.1)) for $E(x+U)-E(x)$,
and then (3.1) will follow.

\medskip
To prove (2.4) we need a precise expression for $I_1(t+G,G)-I_1(t,G)$.
One way to proceed
is to start from (3.2) and use Atkinson's formula (1.3)--(1.7). In the course of
the proof, various expressions will be simplified by Taylor's formula, and one has
to use the well known integral (see e.g., the Appendix of [3])
$$
\int_{-\infty}^\infty \exp\left(Ax-Bx^2\right)\d x = \sqrt{\pi\over B}
\exp\left({A^2\over4B}\right)\qquad(\R B >0).\leqno(3.4)
$$
However, it seems more expedient to proceed in the following way. We start
from Y. Motohashi's formula [17, eq. (4.1.16)], namely
$$\leqno(3.5)$$
\begin{eqnarray*}
{Z}_1(g) &=& \int_{-\infty}^\infty\left[ \R \left\{{\G'\over\G}(\hf+it)\right\}
+2\gamma - \log(2\pi)\right]g(t)\d t + 2\pi\R\{g(\hf i)\}\\
&+ &4\sum_{n=1}^\infty d(n)\int_0^\infty (y(y+1))^{-1/2}g_c(\log(1+1/y))
\cos(2\pi ny)\d y,
\end{eqnarray*}
where, for $k\in\NN$ fixed,
$$
{Z}_k(g) = \int_{-\infty}^\infty |\zt|^{2k}g(t)\d t,\quad
g_c(x) = \int_{-\infty}^\infty g(t)\cos(xt)\d t.\leqno(3.6)
$$
Here the function $g(r)$ is real for $r$ real,
and there exist a large positive constant
$A$ such that $g(r)$ is regular and $g(r) \ll (|r|+1)^{-A}$ for $|\Im m\, r|\le A$.
The choice in our case is
$$
g(t) := {1\over\sqrt{\pi}}\exp\left(-((T-t)/G)^2\right),\quad g_c(x) =
G\exp\left(-{\txt{1\over4}}(Gx)^2\right)\cos(Tx),
$$
and one verifies without difficulty that the above assumptions on $g$ are
satisfied. In this case ${Z}_1(g)$ becomes our $I_1(T,G)$ (see (2.4)).
To evaluate the integral on the right-hand side of (3.5) we use (1.17) with $k=1$.
Thus when we form
the difference $I_1(t+G,G)-I_1(t,G)$ in this way, the integral
on the right-hand side
of (3.5) produces the terms $O(G^2/T) + O(1) = O(1)$, since $G^2/T \ll 1$.
The second integral on the right-hand side of (3.5) is evaluated by the
saddle-point method (see e.g., [3, Chapter 2]). A similar
analysis was made in [8] by the author, and an explicit formula
for $I_1(T,G)$ is also to be found in Y. Motohashi [17, eq. (5.5.1)]. As these
results are either less accurate, or hold in a more restrictive range of $G$ than
what we require for the proof of Theorem 1, a more detailed analysis
is in order. A convenient result to
use is [3, Theorem 2.2 and Lemma 15.1] (due originally to Atkinson [1])
for the evaluation of exponential integrals $\int_a^b\f(x)\exp(2\pi iF(x))\d x$,
where $\f$ and $F$ are suitable smooth, real-valued functions.
 In the latter only the exponential
factor $\exp(-{1\over4}G^2\log(1+1/y))$ is missing. In the notation
of [1] and [3] we have that the saddle point $x_0$ (root of $F'(x)=0$) satisfies
$$
x_0 = U - {1\over2} = \left({T\over2\pi n} +
{1\over4}\right)^{1/2} - {1\over2},
$$
and the presence of the above exponential factor makes it possible to
truncate the series in (3.5) at $n = TG^{-2}\log T$ with a negligible error.
Furthermore, in the remaining range for $n$ we have (in the notation of [3])
$$
\Phi_0\mu_0F_0^{-3/2} \ll (nT)^{-3/4},
$$
which makes a total contribution of $O(1)$, as does error term
integral in Theorem 2.2 of [3]. The error terms with
$\Phi(a),\,\Phi(b)$ vanish for $a \to 0+,\,b \to +\infty\,$.
In this way we obtain a formula, which naturally has a resemblance to
Atkinson's formula (compare it also to [8, eq. (19)]). This is
$$\leqno(3.7)$$
\begin{eqnarray*}
&&I_1(t+G,G)-I_1(t,G) = O(1) +\\& +&
\sqrt{2}G\sum_{n\le TG^{-2}L}(-1)^nd(n)n^{-1/2}\Bigl\{u(t+G,n)H(t+G,n)
- u(t,n)H(t,n)\Bigr\},
\end{eqnarray*}
where
$$
u(t,n) := \left\{\left({t\over2\pi n} + {1\over4}\right)^{1/2} - {1\over2}
\right\}^{-1/2}\qquad\Bigl(t\asymp T, 1 \le n \le TG^{-2}L\Bigr),
$$
and (in the notation of (1.6))
$$
H(T,n) := \exp\left(-G^2\left({\rm {ar\,sinh}}\,\sqrt{\pi n\over2T}\,\right)^2
\right)\sin f(T,n) \quad\Bigl(t\asymp T, 1 \le n \le TG^{-2}L\Bigr).
$$
Now we square (3.7) and integrate over $T \le t\le 2T$,
using the Cauchy-Schwarz inequality and (3.1) to obtain
$$
\int_T^{2T}\Bigl(I_1(t+G,G)-I_1(t,G)\Bigr)^2\d t \;=\;
 S + O(TG^{1/2}L^{3/2}),\leqno(3.8)
$$
where we set
$$\leqno(3.9)
$$
\begin{eqnarray*}
S &:=& 2G^2\int_T^{2T}\Biggl\{\sum_{n\le TG^{-2}L}(-1)^n
\frac{d(n)}{n^{1/2}}\Bigl[u(t+G,n)H(t+G,n)- \\
&-& u(t,n)H(t,n)\Bigr]\Biggr\}^2\d t.
\end{eqnarray*}
Squaring out the sum over $n$ in (3.9) it follows that
$$\leqno(3.10)
$$
\begin{eqnarray*}
S &=& 2G^2\int_T^{2T}\sum_{n\le TG^{-2}L}{d^2(n)\over n}
\Bigl[u(t+G,n)H(t+G,n) - u(t,n)H(t,n)\Bigr]^2\d t + \\
&+& 2G^2\int\limits_T^{2T}\sum_{m\ne n\le TG^{-2}L}(-1)^{m+n}{d(m)d(n)\over\sqrt{mn}}
\Bigl[u(t+G,m)H(t+G,m) \,-  \\ & -&u(t,m)H(t,m)\Bigr]
\Bigl[u(t+G,n)H(t+G,n) - u(t,n)H(t,n)\Bigr]\d t.
\end{eqnarray*}
The main term in Theorem 1 comes from the first sum in (3.10) (the diagonal terms),
while the sum over the non-diagonal terms $m\ne n$ will contribute
to the error term. To see this  note first that
the functions $u(t,n)\; (\asymp (n/t)^{1/4}$ in our range) and
$$
\exp\left(-G^2\left({\rm {ar\,sinh}}\,\sqrt{\pi n\over2T}\,\right)^2
\right)
$$
are monotonic functions of $t$ when $t\in [T,\,2T]\,$, and moreover, since
$$
{\partial f(t,n)\over \partial t} \;=\; 2{\rm {ar\,sinh}}\,\sqrt{\pi n\over2t}\,,
$$
it follows that, for $U,V = 0$ or $G$,
$$
{\partial [f(t+U,m)\pm f(t+V,n)]\over \partial t} \;\asymp\; {|\sqrt{m} \pm
\sqrt{n}|\over\sqrt{T}}\qquad(m\ne n). \leqno(3.11)
$$
In the sum over $m\ne n$ we can assume, by symmetry, that $n < m\le 2n$ or $m>2n$.
Hence by the first derivative test (see e.g., Lemma 2.1 of [3]) and (3.11) we have
that the sum in question is
\begin{eqnarray*}
&\ll& G^2\sum_{1\le n< m\le TG^{-2}L}{d(m)d(n)\over(mn)^{1/4}}\cdot{1\over\sqrt{m}
-\sqrt{n}}\\
&\ll& G^2\sum_{n\le TG^{-2}L}d(n)\sum_{n<m\le2n}{d(m)\over m-n} +\\
&+& G^2\sum_{m\le TG^{-2}L}{d(m)\over m^{3/4}}\sum_{n<2m}{d(n)\over n^{1/4}}\\
&
\ll& G^2T^\e\Bigl(\sum_{n\le TG^{-2}L}1 +
\sum_{m\le TG^{-2}L}d(m)L\Bigr) \ll_\e T^{1+\e}.
\end{eqnarray*}
Note that, by the mean value theorem,
$$
u(t+G,n) - u(t,G) = O(Gn^{1/4}T^{-5/4})\qquad(t\asymp T, 1 \le n \le TG^{-2}L).
$$
Hence we obtain, by trivial estimation,
$$\leqno(3.12)
$$
\begin{eqnarray*}
S &= &2G^2\int_T^{2T}\sum_{n\le TG^{-2}L}{d^2(n)\over n}u^2(t,n)
\Bigl(H(t+G,n)-H(t,n)\Bigr)^2\d t \\
&+ &O\Biggl(G^2\int_T^{2T}\sum_{n\le TG^{-2}L}{d^2(n)\over n}
\cdot Gn^{1/2}T^{-3/2}\d t\Biggr)
+ O_\e(T^{1+\e})\\
&=& 2G^2\int_T^{2T}\sum_{n\le TG^{-2}L}{d^2(n)\over n}u^2(t,n)
\Bigl(H(t+G,n)-H(t,n)\Bigr)^2\d t  + O_\e(T^{1+\e}),
\end{eqnarray*}
since $G^2 \ll T$. Now note that
$$
u^2(t,n) = \left({t\over2\pi n}\right)^{-1/2} + O\left({n\over T}\right)
\qquad(T\le t\le 2T),
$$
hence the error term above makes a contribution to $S$ which is
$$
\ll G^2\sum_{n\le TG^{-2}L}d^2(n) \ll TL^4.
$$
Similarly, replacing $t+G$ by $t$ in the exponential factor in $H(t+G,n)$, we make
a total error which is $\ll_\e T^{1+\e}$. Therefore (3.8) and (3.12) give
\pagebreak
$$\leqno(3.13)
$$
\begin{eqnarray*}
&&\int_T^{2T}\Bigl(I_1(t+G,G)-I_1(t,G)\Bigr)^2\d t \;=\;  O(TG^{1/2}L^{3/2})\\
&+& 2\sqrt{2\pi}G^2\sum_{n\le TG^{-2}L}{d^2(n)\over\sqrt{n}}\int_T^{2T}
t^{-1/2}\exp\left(-2G^2\left({\rm {ar\,sinh}}\,
\sqrt{\pi n\over2t}\,\right)^2\right)\times
\\
&&\times\Bigl(\sin f(t+G,n) - \sin f(t,n)\Bigr)^2\d t
\quad\qquad(T^\e\le G\ll\sqrt {T}).
\end{eqnarray*}
In a similar vein we simplify  (3.13), by using
$$
{\rm {ar\,sinh}}\,z = z + O(|z|^3)\qquad(|z| < 1).
$$
Thus we may replace the exponential factor in (3.13) by $\exp(-\pi nG^2/t)$, making
an error which is absorbed by the $O$--term in (3.13). Next we use the identity
$$
\Bigl(\sin\a - \sin \b\Bigr)^2 = (2+2\cos(\a+\b))\sin^2\hf(\a-\b)
\quad( \a = f(t+G,n),\;\b = f(t,n)).
$$
Analogously to the treatment of the sum in (3.10) with $m\ne n$,
we use the first derivative test
to show that the contribution of the terms with $\cos(\a+\b)$ is
$O(TG^{1/2}L^{3/2})$. Therefore (3.13) reduces to
$$\leqno(3.14)
$$
\begin{eqnarray*}
&&\int_T^{2T}\Bigl(I_1(t+G,G)-I_1(t,G)\Bigr)^2\d t \;=\;  O(TG^{1/2}L^{3/2})\\
&
+& 4\sqrt{2\pi}G^2\sum_{n\le TG^{-2}L}{d^2(n)\over\sqrt{n}}\int_T^{2T}
t^{-1/2}\exp\left(-{\pi nG^2\over t}\right)
\sin^2\left(\sqrt{\pi n\over 2t}G\right)\d t,
\end{eqnarray*}
since
$$
\sin^2\hf(\a-\b) = \sin^2\left(\sqrt{\pi n\over 2t}G\right) +
O\left((G^2n^{1/2} + Gn^{3/2})T^{-3/2}\right).
$$
In the integral on the right-hand side of (3.14) we make the change of
variable
$$
\sqrt{\pi n\over 2t}G = y, \; t = {\pi nG^2\over 2y^2},\;
\d t = - {\pi nG^2\over y^3}\d y.
$$
The main term on the right-hand side of (3.14) becomes then
\pagebreak
$$\leqno(3.15)
$$
\begin{eqnarray*}
&&8\pi G^3\sum_{n\le TG^{-2}L}d^2(n)\int_{\sqrt{\pi n\over 4T}G}
^{\sqrt{\pi n\over 2T}G}\left({\sin y\over y}\right)^2{\rm e}^{-2y^2}\d y\\
&
=& 8\pi G^3\int_{\sqrt{\pi \over 4T}G}^{\sqrt{\pi L\over2}}
\sum_{\max(1,{2Ty^2\over\pi G^2})\le n\le \min(TG^{-2}L,{4Ty^2\over\pi G^2})}d^2(n)
\cdot \left({\sin y\over y}\right)^2{\rm e}^{-2y^2}\d y\\
&
=& 8\pi G^3\int_{\sqrt{\pi\over2T}G}^{{1\over2}\sqrt{\pi L}}
\sum_{{Ty^2\over\pi G^2}\le n \le {2Ty^2\over\pi G^2}}d^2(n)
\left({\sin y\over y}\right)^2{\rm e}^{-2y^2}\d y + O(T^{1/2}G^2L^4).
\end{eqnarray*}
At this point we invoke (see [8] for a proof, and [16] for a slightly sharper result)
the asymptotic formula
$$
\sum_{n\le x}d^2(n) = x\sum_{j=0}^3d_j\log^jx + O_\e(x^{1/2+\e})\qquad(d_3 =
1/(2\pi^2)).\leqno(3.16)
$$
By using (3.16) it follows ($a_j, b_j,b'_j,d_j$ denote constants which may be explicitly
evaluated) that the last main term in (3.15)  equals
\begin{eqnarray*}&{}&
8\pi G^3\int_{\sqrt{\pi\over2T}G}^{{1\over2}\sqrt{\pi L}}
\Biggl\{{Ty^2\over G^2}\sum_{j=0}^3b_j\log^j\Bigl({Ty^2\over G^2}\Bigr)\\
&
+& O_\e\left({T^{1/2+\e}y\over G}\right)\Biggr\}
\left({\sin y\over y}\right)^2{\rm e}^{-2y^2}\d y +O(T^{1/2}G^2L^4)\\
&
=& 8\pi TG\int_0^\infty \sin^2y\left(\sum_{j=0}^3b'_j\log^j
\Bigl({\sqrt{T}y\over G}\Bigr)\right)
{\rm e}^{-2y^2}\d y +O_\e(T^{1/2+\e}G^2)\\
&
=& TG\sum_{j=0}^3a_j\log^j\Bigl({\sqrt{T}\over G}\Bigr) + O_\e(T^{1/2+\e}G^2).
\end{eqnarray*}
Coupled with (3.14)--(3.15) this proves Theorem 1 with
$$
a_3 = 8b'_3\pi\int_0^ \infty {\rm e}^{-2y^2}\sin^2y \d y > 0.
$$
Namely by using (3.4) we have
\begin{eqnarray*}
\int_0^\infty {\rm e}^{-2y^2}\sin^2y\d y &=&
\txt{1\over4}\int_{-\infty}^\infty {\rm e}^{-2y^2}(1-\cos 2y)\d y \\
&
=& \txt{1\over4}\R \Biggl\{\int_{-\infty}^\infty {\rm e}^{-2y^2}
(1- {\rm e}^{2iy})\d y \Biggr\}\\&
=&
{\sqrt{\pi}\over4\sqrt{2}}\left(1- {1\over\sqrt{\rm e}}\right),
\end{eqnarray*}
and the other constants $a_j$ in (2.4) can be also explicitly evaluated.
This finishes then the proof of Theorem 1.

\section{A generalization of the mean square result}

In Theorem 1  we considered the mean square of $I_1(t+G,G) - I_1(t,G)$ (see (2.4)),
over the ``long" interval $[T,\,2T]\,$. However, already M. Jutila [13]
(see (1.8) and (1.9)) considered the mean square of $\D(x+U)-\D(x)$
and $E(t+U)-E(t)$ over the ``short" interval $[T,\,T+H]$.
Therefore it seems also natural to consider the mean square of
$I_1(t+U,G) - I_1(t,G)$ over the short interval $[T,\,T+H]$ for suitable $U = U(T)$.
It turns out that this problem is more complicated, because of the
presence of two parameters $U$ and $G$, and not only $U$ as
in (1.8) and (1.9). Our assumption will be henceforth that
$$
T^\e \le U = U(T) \le G = G(T) \le \hf\sqrt{T},\;
T^\e \le H =H(T) \le T, \;HU \gg T^{1+\e}.
\leqno(4.1)
$$
The method of proof will be analogous to the proof of Theorem 1,
only it will be technically  more involved, and the final
result will not have such a nice shape as (2.5). We shall thus
only sketch the salient points of the evaluation of
$$
J(T) = J(T;G,H,U) := \int_T^{T+H}{\Bigl((I_1(t+U,G) - I_1(t,G)\Bigr)}^2\d t,
\leqno(4.2)
$$
subject to the condition (4.1), without formulating a theorem.

\medskip
To obtain an upper bound for $J(T)$ we recall (3.2), which gives now
$$\leqno(4.3)
$$
\begin{eqnarray*}
&&I_1(t+U,G)-I_1(t,G)= {2\over\sqrt{\pi}G^2}\int\limits_{-GL}^{GL}u{\rm e}^{-(u/G)^2}
\Bigl(E(t+u+G)-E(t+u)\Bigr)\d u  \\
& +& O\Bigl({UG\over T}\Bigr)
= {2\over\sqrt{\pi}G^2}\int\limits_{-GL}^{GL}u{\rm e}^{-(u/G)^2}
\Bigl(E(t+u+G)-E(t+u)\Bigr)\d u + O(1),
\end{eqnarray*}
since $U \le G \ll \sqrt{T}$. First we square (4.3), integrate
over $[T,\,T+H]$ and use (1.9) for $E(t)$ to obtain
$$
J(T) = \int_T^{T+H}{\Bigl((I_1(t+U,G) - I_1(t,G)\Bigr)}^2\d t
\ll_\e HU\log^3{\sqrt{T}\over U}.\leqno(4.4)
$$
Now we use (3.7), square, integrate over $[T,\,T+H]$ and use (4.4). It follows
that
$$
J(T) = \int_T^{T+H}{\Bigl(I_1(t+U,G) - I_1(t,G)\Bigr)}^2\d t
= {\bf S} +  O(H\sqrt{U}L^{3/2}),\leqno(4.5)
$$
where, similarly to (3.9), now we shall have
\begin{eqnarray*}
{\bf S} \!\!&:=& 2G^2\int_T^{T+H}\Biggl\{\sum_{n\le TG^{-2}L}(-1)^nd(n)n^{-1/2}
\Bigl[u(t+U,n)H(t+U,n)\\
& \;-& u(t,n)H(t,n)\Bigr]\Biggr\}^2\d t.
\end{eqnarray*}
Proceeding as in the proof of (3.10)--(3.13) we shall obtain
$$\leqno(4.6)$$
\begin{eqnarray*}
{\bf S}\!\!&=&
2\sqrt{2\pi}G^2\sum_{n\le TG^{-2}L}{d^2(n)\over\sqrt{n}}\int_T^{T+H}
t^{-1/2}\exp\left(-2G^2\left({\rm {ar\,sinh}}\,
\sqrt{\pi n\over2t}\,\right)^2\right)\times
\\
&&\times\Bigl(\sin f(t+U,n) - \sin f(t,n)\Bigr)^2\d t
 +  O(H\sqrt{U}L^{3/2}) +  O_\e(T^{1+\e}).
\end{eqnarray*}
Using again the identity
$$
\Bigl(\sin\a - \sin \b\Bigr)^2 = (2+2\cos(\a+\b))\sin^2\hf(\a-\b)
\quad( \a = f(t+U,n),\;\b = f(t,n))
$$
and simplifying the exponential factor in (4.6) by Taylor's formula, we have that
$$
{\bf S} = 4\sqrt{2\pi}G^2\sum_{n\le TG^{-2}L}{d^2(n)\over\sqrt{n}}
\int\limits_T^{T+H}t^{-1/2}\exp\left(-{\pi nG^2\over t}\right)
\sin^2\hf(\a-\b)\d t + O_\e(R),\leqno(4.7)
$$
where
$$
R := H\sqrt{U}L^{3/2} +  T^{1+\e}.
$$
Note that
$$
\sin^2\hf(\a-\b) = \sin^2\left({\sqrt{\pi n\over2t}}\,U\right)
+O\left((U^2n^{1/2} + Un^{3/2})T^{-3/2}\right)
$$
in the relevant range, and the total contribution of the $O$--terms above
will be $O_\e(T^\e H)$. In the integral in (4.7) we make the change of variable,
similarly as was done on the right-hand side of (3.14),
$$
\sqrt{\pi n\over t}G = y, \; t = {\pi nG^2\over y^2},\;
\d t = - {2\pi nG^2\over y^3}\d y.
$$
The main term in (4.7) becomes, after changing the order of integration
and summation,  $O_\e(T^\e R)$ plus
$$
8\pi\sqrt{2}G^3\int_{G\sqrt{\pi/T}}^{\sqrt{\pi LT/(T+H)}}
\sum_{{Ty^2\over\pi G^2}\le n\le {(T+H)y^2\over\pi G^2}}
d^2(n)\sin^2\left({U\over G\sqrt{2}}y\right)y^{-2}{\rm e}^{-y^2}\d y.
$$
For the sum over $n$ we use again the asymptotic formula (3.16).
We obtain that
\pagebreak
$$
\leqno(4.8)
$$
\begin{eqnarray*}{\bf S}\!\!&=\!&
8\pi\sqrt{2}G^3\int_{G\sqrt{\pi/T}}^{\sqrt{\pi LT/(T+H)}}
\left\{xP_3(\log x)\Biggl|_{x=Ty^2/(\pi G^2)}^{x=(T+H)y^2/(\pi G^2)}
+ O\left({T^{1/2+\e}y\over G^2}\right)\right\}\times\\
&&
\times \sin^2\left({U\over G\sqrt{2}}y\right)y^{-2}{\rm e}^{-y^2}\d y
+ O_\e(T^{1+\e}) + O(H\sqrt{U}L^{3/2}),
\end{eqnarray*}
where (cf. (3.16)) $P_3(z) = \sum_{j=0}^3d_jz^j$. The main term in (4.8)
equals
$$\leqno(4.9)
$$
\begin{eqnarray*}
&{}&
8\sqrt{2}G\int\limits_{G\sqrt{\pi/T}}^{\sqrt{\pi LT/(T+H)}}
xP_3\left(\log\left({xy^2\over \pi G^2}\right)\right)\Biggl|_{x=T}^{x=T+H}
\sin^2\left({U\over G\sqrt{2}}y\right){\rm e}^{-y^2}\d y\\
&
=& 8\sqrt{2}G\int_{0}^{\infty}
xP_3\left(\log\left({xy^2\over \pi G^2}\right)\right)\Biggl|_{x=T}^{x=T+H}
\sin^2\left({U\over G\sqrt{2}}y\right){\rm e}^{-y^2}\d y\\
&
+& O_\e(T^{\e-3/2}HU^2G^2).
\end{eqnarray*}
In view of (4.1) the last $O$--term is $\ll_\e T^{1/2+\e}U^2$, so that (4.9)
gives
\begin{eqnarray*}
{\bf S}\!\!\!&=&\!\!\!
8\sqrt{2}G\int_{0}^{\infty}
xP_3\left(\log\left({xy^2\over \pi G^2}\right)\right)\Biggl|_{x=T}^{x=T+H}
\sin^2\left({U\over G\sqrt{2}}y\right){\rm e}^{-y^2}\d y\\
&
+& O_\e(T^{1+\e})+ O_\e(T^{1/2+\e}U^2)+ O(H\sqrt{U}L^3).
\end{eqnarray*}
This can be rewritten as
$$
\leqno(4.10)
$$
\begin{eqnarray*}
{\bf S}\!\!\!&=&\!\!\!
G\int_{0}^{\infty}\left\{x\sum_{k=0}^3A_k(y)\log^k\left({\sqrt{x}\over G}\right)
\right\}
\Biggl|_{x=T}^{x=T+H}
\sin^2\left({U\over G\sqrt{2}}y\right){\rm e}^{-y^2}\d y\\
&+& O_\e(T^{1+\e})+ O_\e(T^{1/2+\e}U^2)+ O(H\sqrt{U}L^3),
\end{eqnarray*}
where
$$
A_k(y) = A_k(y;U,G) := \sum_{j=0}^3 b_{j,k}\left(\log {y^2\over\pi}\right)^j
$$
with computable coefficients $b_{j,k} \,(= b_{j,k}(U,G))$.
This shows that the main term in  $\bf S$ has the form
$$
Gx\sum_{k=0}^3D_k\log^k\left({\sqrt{x}\over G}\right)\Biggl|_{x=T}^{x=T+H}
$$
with computable coefficients $D_k = D_k(U,G)$.
In particular, by using (3.4) we see that $D_3$ is a multiple of
\begin{eqnarray*}
&{}&\int_{0}^{\infty}\sin^2\left({U\over G\sqrt{2}}y\right){\rm e}^{-y^2}\d y
= {1\over4}\int_{-\infty}^{\infty}
\left(1- \cos\Bigl({\sqrt{2}U\over G}y\Bigr)\right){\rm e}^{-y^2}\d y
\\&
= &{1\over4}\R\left\{\int_{-\infty}^{\infty}
\Bigl(1 - \exp\Bigl(i{\sqrt{2}U\over G}y\Bigr)\Bigr)
{\rm e}^{-y^2}\d y\right\} = {\sqrt{\pi}\over4}\left(1 - {\rm e}^{-U^2/(2G^2)}
\right),
\end{eqnarray*}
that is, $D_3$ is an expression depending only on $T$.
For $U =o(G)$ we have $D_3 = (C+o(1))U^2G^{-2}$
 (with $C>0$ and $T\to\infty$). This shows that, if the parameters $U, G, H$
are suitably chosen as functions of $T$, then (4.5)--(4.10) give
$$
\int_T^{T+H}{\Bigl(I_1(t+U,G)-I_1(t,G)\Bigr)}^2 \d t\;\asymp\;
{HU^2\over G}\log^3\Biggl({\sqrt{T}\over G}\Biggr)\qquad(T\to\infty),\leqno(4.11)
$$
which is more precise than (4.4). It is clear that, in that case,
(4.11) can be in fact replaced by
a true asymptotic formula (for example, $U = T^{1/3}, G = T^{4/9}, H = T^{8/9}$
is such a choice) for $J(T)$. Such a formula can be written down explicitly,
although its form will be unwieldy, and because of this it is not
formulated as a theorem.

\section{The proof of Theorem 2}

We assume that the hypotheses of Theorem 2 hold, namely that
$$
T^\e \le U \le GT^{-\e}\ll T^{1/2-\e}, \quad U = U(T),\quad G=G(T).\leqno(5.1)
$$
We start from (1.14)--(1.16) and use (1.17) to deduce that, for $T\le t\le2T,\,
Y_0 = Y_0(t;\k_j)$,
$$
\leqno(5.2)
$$
\begin{eqnarray*}&{}&
I_2(t+U,G)-I_2(t,G) \sim O(UG^{-1/2}L^C)+ O(1) \;+\\&
+ &{G\pi\over\sqrt{2t}}\sum_{\kappa_j\le TG^{-1}L}\a_j H_j^3(\hf)
 \kappa_j^{-1/2}{\rm e}^{-{1\over4}G^2\log^2(1+Y_0)}
 \bigl(f_j(t+U)-f_j(t)\bigr),
 \end{eqnarray*}
where
$$
f_j(T) :=\sin\Bigl(\kappa_j \log
{\kappa_j\over 4{\rm e}T} +c_3\k_j^3T^{-2}+ \cdots + c_N\k_j^NT^{1-N}\Bigr).
\leqno(5.3)
$$
Here we used the bounds
\begin{eqnarray*}&&
(t+U)^{-1/2} - t^{-1/2} \;\ll\; UT^{-3/2},\\&&
{\rm e}^{-{1\over4}G^2\log^2\bigl(1+Y_0(t+U;\k_j)\bigr)}-
{\rm e}^{-{1\over4}G^2\log^2\bigl(1+Y_0(t;\k_j)\bigr)}
\ll UT^{-1}L^2,
\end{eqnarray*}
which follows from Taylor's formula, $\k_j \ll TG^{-1}L$,
and (see Y. Motohashi [17])
$$
\sum_{K<\k_j\le2K}\a_j H_j^3(\hf) \;\ll\; K^2\log^CK.\leqno(5.4)
$$
Taylor's formula yields, for a fixed integer $L\ge 1$ and some
$\theta$ satisfying $|\theta| \le 1$,
$$
f_j(t+U)-f_j(t) = \sum_{\ell=1}^L{U^\ell\over\ell !}f_j^{(\ell)}(t)
+ O\left({U^{L+1}\over (L+1)!}\Bigl|f_j^{(L+1)}(t+\theta U)\Bigr|\right).
\leqno(5.5)
$$
Observe that
\begin{eqnarray*}
f_j'(t) &=& \cos\Bigl(\kappa_j \log
{\kappa_j\over 4{\rm e}t} +\cdots\Bigr)\left(-{\k_j\over t} -
2c_3\k_j^3t^{-3}-\cdots\right), \\
f_j''(t) &=& -\sin\Bigl(\kappa_j \log
{\kappa_j\over 4{\rm e}t} +\cdots\Bigr)\left(-{\k_j\over t} -
2c_3\k_j^3t^{-3}-\cdots\right)^2 + \\&
+& \cos\Bigl(\kappa_j \log
{\kappa_j\over 4{\rm e}t} +\cdots\Bigr)\left({\k_j\over t^2} + 6c_3\k_j^3t^{-4}
+\cdots\right),
\end{eqnarray*}
and so on. Since $U \ll GT^{-\e}$, this means that for
$L\;(= L(\e))$ sufficiently large
the last term in (5.5) makes, by trivial estimation, a negligible contribution
(i.e., $\ll 1$). Each time the derivative is decreased by a factor which is
$$
\ll U\k_j T^{-1} \ll UTG^{-1}LT^{-1} \ll UG^{-1}L \ll_\e T^{-\e/2}.
$$
This means that in (5.5) the term $\ell=1$, namely $Uf_j'(t)$ will make
the largest contribution. This contribution is, on squaring (5.2)
and integrating,
$$
\ll {G^2L\over T}\max_{K\ll TG^{-1}L}\int_{T/2}^{5T/2}\f(t)
\Bigl|\sum(K)\Bigr|^2\d t.
$$
Here $\f(t)\;(\ge0)$ is a smooth function supported on $[T/2,\,5T/2]$ such
that $\f(t) =1$ when $T\le t \le2T$ and $\f^{(r)}(t)\ll_r T^{-r}$ for
$r = 0,1,2,\ldots$, and
\begin{eqnarray*}
\sum(K) \!\!\!&:=&\!\!\! U\sum_{K<\k_j\le K'\le2K}\a_j H_j^3(\hf)
 \k_j^{-1/2}{\rm e}^{-{1\over4}G^2\log^2(1+Y_0(t;\k_j))}\\&\times&
\!\!\!\Bigl({\k_j\over t} + 2c_3{\k_j^3\over t^{3}}
+\cdots\Bigr)\cos\Bigl(\kappa_j \log
{\kappa_j\over 4{\rm e}t} +\cdots\Bigr).
\end{eqnarray*}
When $\sum(K)$ is squared, we shall obtain a double sum over $K <\k_j,\k_\ell\le K'$,
with the exponential factors (both of which are estimated analogously)
$$
\exp\Bigl(if_j(t)-if_\ell(t)\Bigr),\quad \exp\Bigl(if_\ell(t)-if_j(t)\Bigr),
$$
in view of (5.3). The first one yields the integral
$$
I:= \int_{T/2}^{5T/2}{\rm e}^{-{1\over4}G^2\log^2(1+Y_0(t;\k_j))}
{\rm e}^{-{1\over4}G^2\log^2(1+Y_0(t;\k_\ell))}
F(t;\k_j,\k_\ell)t^{i\k_\ell-i\k_j}\d t,
$$
where for brevity we set
$$
F(t;\k_j,\k_\ell):= \f(t)
 \exp \left\{i\Bigl(c_3(\k_j^3-\k_\ell^3)t^{-2} +
\cdots + c_N(\k_j^N-\k_\ell^N)t^{1-N}\Bigr)\right\}.
$$
Integration by parts shows that
$$
I = -  \int_{T/2}^{5T/2}
{\left\{{\rm e}^{-\ldots}{\rm e}^{-\ldots}
F(t;\k_j,\k_\ell)\right\}}'{t^{i\k_\ell-i\k_j+1}
\over i\k_\ell- i\k_j+1}\d t,
$$
and
$$
{\left\{{\rm e}^{-\ldots}{\rm e}^{-\ldots}
F(t;\k_j,\k_\ell)\right\}}'
 \;\ll\; {1\over T} + {|\k_j-\k_\ell|K^2\over T^3}.
$$
Therefore if integration by parts is performed a sufficiently large number of times
(depending on $\e$),
then the contribution of $\k_j,\k_\ell$ which satisfy $|\k_j-\k_\ell|\ge T^\e$
will be negligibly small, since the integrand is decreased each time by a factor
which is, for $|\k_j-\k_\ell|\ge T^\e$,
$$
\ll\; {T\over|\k_j- \k_\ell+1|}\left({1\over T} + {|\k_j-\k_\ell|K^2\over T^3}\right)
\;\ll\; T^{-\e},
$$
and the exponential factor remains the same.
To complete the proof of Theorem 2 we use the bound, proved by the author in [5],
$$
\sum_{K\le \k_j\le K+1}\a_jH_j^3(\hf) \;\ll_\e\; K^{1+\e}.\leqno(5.6)
$$
From (5.6) and the preceding discussion it follows that
\begin{eqnarray*}
&&\int_{T/2}^{5T/2}\f(t)
\Bigl|\sum(K)\Bigr|^2\d t\\&
\ll& {U^2K^2\over T^2}T\sum_{K<\k_j\le K'}\a_j\k_j^{-1/2}H_j^3(\hf)
\sum_{|\k_j-\k_\ell|\le T^\e}\a_\ell \k_\ell^{-1/2} H_\ell^3(\hf)\\&
\ll_\e& {U^2K^2\over T}KT^\e K^2K^{-1} \ll_\e U^2T^{3+\e}G^{-4},
\end{eqnarray*}
which gives for the integral in (2.6) the bound
$$
TU^2G^{-1}L^C + TL^8 + T^{2+\e}\left({U\over G}\right)^2.
$$
However, it is clear that in our range
$$
TU^2G^{-1}L^C \;\ll\; T^{2+\e}(U/G)^2
$$
holds. Moreover (5.1) implies that
$$
T^{2+\e}\left({U\over G}\right)^2 \;\gg\; T^{2+4\e}G^{-2} \;\ge\; T^{1+4\e},
$$
so that the bound in (2.6) follows. Thus Theorem 2 is proved, but the true order
of the integral in (2.6) is elusive. Namely in the course of the proof we estimated,
by the use of (5.6), trivially an exponential sum with $\a_jH_j^3(\hf)$. This
certainly led to some loss, and for the discussion of bounds for
exponential sums with
Hecke series, of the type needed above,
 the reader is referred to the author's recent work [9].

\end{document}